\documentclass[11pt]{amsart}

\usepackage{amsmath}
\usepackage{amssymb}
\usepackage{mathrsfs}

\newtheorem{theorem}{Theorem}[section]
\newtheorem{claim}[theorem]{Claim}

\newtheorem{proposition}[theorem]{Proposition}
\newtheorem{corollary}[theorem]{Corollary}

\theoremstyle{definition}

\newtheorem{question}[theorem]{Question}

\theoremstyle{remark}

\newcount\skewfactor
\def\mathunderaccent#1#2 {\let\theaccent#1\skewfactor#2
\mathpalette\putaccentunder}
\def\putaccentunder#1#2{\oalign{$#1#2$\crcr\hidewidth
\vbox to.2ex{\hbox{$#1\skew\skewfactor\theaccent{}$}\vss}\hidewidth}}
\def\name{\mathunderaccent\tilde-3 }

\def\smallbox#1{\leavevmode\thinspace\hbox{\vrule\vtop{\vbox
   {\hrule\kern1pt\hbox{\vphantom{\tt/}\thinspace{\tt#1}\thinspace}}
   \kern1pt\hrule}\vrule}\thinspace}

\newcommand{\cf}{{\rm cf}}

\def\qedref#1{$\qed_{\reforiginal{#1}}$}

\title{Unbalanced polarized relations}
\author{Shimon Garti}
\address{Einstein Institute of Mathematics,
 The Hebrew University of Jerusalem,
 Jerusalem 91904, Israel}
\email{shimon.garty@mail.huji.ac.il}


\begin{document}
\let\labeloriginal\label
\let\reforiginal\ref
\def\ref#1{\reforiginal{#1}}
\def\label#1{\labeloriginal{#1}}

\begin{abstract}
We prove the consistency of the relation $\binom{\mu^+}{\mu}\nrightarrow\binom{\mu^+\ \omega_1}{\mu\quad \mu}$ where $\mu$ is a strong limit singular cardinal of countable cofinality.
This result can be forced at $\mu=\aleph_\omega$.
\end{abstract}

\maketitle

\newpage

\section{Introduction}

The polarized partition relation $\binom{\alpha}{\beta}\rightarrow \binom{\gamma_0\ \gamma_1}{\delta_0\ \delta_1}$ says that for every coloring $c:\alpha\times\beta\rightarrow 2$ there are $A\subseteq\alpha,B\subseteq\beta$ and $i\in\{0,1\}$ such that ${\rm otp}(A)=\gamma_i,{\rm otp}(B)=\delta_i$ and $c\upharpoonright(A\times B)$ is constantly $i$.
If $(\gamma_0,\delta_0)\neq(\gamma_1,\delta_1)$ then we shall say that the relation is \emph{unbalanced}.
Lest $\gamma_0=\gamma_1=\gamma$ and $\delta_0=\delta_1=\delta$ we write $\binom{\alpha}{\beta}\rightarrow\binom{\gamma}{\delta}$ and then we shall say that the relation is \emph{balanced}.

An old problem raised by Erd\H{o}s, Hajnal and Rado in the so-called \emph{Giant triple paper} is whether $\binom{\aleph_{\omega+1}}{\aleph_\omega}\rightarrow\binom{\aleph_{\omega+1}\ \omega_1}{\aleph_\omega\quad \aleph_\omega}$ under \textsf{GCH}.
More generally, the question arises with respect to the relation $\binom{\mu^+}{\mu}\rightarrow\binom{\mu^+\ \omega_1}{\mu\quad \mu}$ where $\mu>\cf(\mu)=\omega$ is strong limit and $2^\mu=\mu^+$, see \cite[Questions 1 and 2]{MR2367118}.
A particular interesting case is the case of an $\omega$-limit of measurable cardinals, see \cite[Question 6]{GMS}.

We suggest a negative answer to these problems.
We shall prove that consistently $\mu>\cf(\mu)=\omega, \mu$ is a strong limit cardinal, $2^\mu=\mu^+$ yet $\binom{\mu^+}{\mu}\nrightarrow\binom{\mu^+\ \omega_1}{\mu\quad \mu}$.
In this result it is possible that $\mu$ be a limit of measurable cardinals or a very small cardinal like $\aleph_\omega$.

Our notation is mostly standard.
We follow \cite{MR795592} with respect to arrows notation.
We adopt the conventions of \cite{MR2768695} with respect to Prikry type forcing, and in particular we use the Jerusalem forcing notation.
We suggest the wonderful monograph \cite{MR3075383} for basic results about polarized partition relations, and the Handbook chapter \cite{MR2768681} for advanced material concerning polarized relations.

I am very much obliged to the referee of this paper.
My original manuscript contained another section with a false proof.
The referee pointed to my mistake and I am sincerely grateful.

\newpage 

\section{Unbalanced relations}

In this section we deal with an unbalanced relation at a strong limit singular cardinal $\mu$ so that $2^\mu=\mu^+$.
Erd\H{o}s, Hajnal and Rado proved in \cite{MR0202613} that if $2^\mu=\mu^+$ then 
$\binom{\mu^+}{\mu}\nrightarrow\binom{\mu^+}{\mu}$ for every infinite cardinal $\mu$.
The positive relation $\binom{\mu^+}{\mu}\rightarrow\binom{\mu^+}{\mu}$ is consistent at strong limit singular cardinals, see \cite{MR2987137}, but this requires $2^\mu>\mu^+$.
It also holds under \textsf{AD}, where the assumption $2^\mu=\mu^+$ becomes irrelevant, see \cite{MR4101445}.
One may wonder, therefore, what is the best positive relation at such cardinals without the assumption $2^\mu>\mu^+$, that is under the assumption $2^\mu=\mu^+$ or just in \textsf{ZFC}.
An elegant result of Shelah from \cite{MR1606515} says that if $\mu$ is a singular cardinal and a limit of measurable cardinals then $\binom{\mu^+}{\mu}\rightarrow\binom{\tau}{\mu}$ for every $\tau\in\mu^+$.
Being a theorem of \textsf{ZFC}, it holds even if $2^\mu=\mu^+$.

These facts lead to the investigation of the intermediate unbalanced relation $\binom{\mu^+}{\mu}\rightarrow\binom{\mu^+\ \tau}{\mu\quad \mu}$.
On the one hand, this relation is weaker than $\binom{\mu^+}{\mu}\rightarrow\binom{\mu^+}{\mu}$.
On the other hand, it is stronger than $\binom{\mu^+}{\mu}\rightarrow\binom{\tau}{\mu}$.
Due to Shelah's result, this intermediate relation is particularly interesting at a singular limit of measurable cardinals.

By simple arguments one can show that $\binom{\mu^+}{\mu}\rightarrow\binom{\mu^+\ n}{\mu\quad \mu}$ for every $n\in\omega$.
However, if $\mu\geq\cf(\mu)>\omega$ then $\binom{\mu^+}{\mu}\nrightarrow\binom{\mu^+\ \omega}{\mu\quad \mu}$ as proved in \cite{MR0202613}.
Consequently, the question is settled for cardinals with uncountable cofinality (both regular and singular).
This is the reason for concentrating on singular cardinals with countable cofinality.

It is quite surprising to find out that for strong limit singular cardinals with countable cofinality under the assumption $2^\mu=\mu^+$ we have $\binom{\mu^+}{\mu}\rightarrow\binom{\mu^+\ \omega}{\mu\quad \mu}$, as proved in \cite{MR0202613}.
Motivated by this peculiar situation in which singular cardinals with countable cofinality demonstrate \emph{a stronger} positive relation, Erd\H{o}s, Hajnal and Rado tried to check how large can the small component in the second color be.
Assuming \textsf{GCH} they proved that $\binom{\mu^+}{\mu}\nrightarrow\binom{\mu^+\ \omega_2}{\mu\quad \mu}$, and later it was shown that \textsf{GCH} can be replaced by $2^\mu=\mu^+$ only, see \cite{GMS}.
The remaining case is, therefore, $\omega_1$ at the second color.
Jones proved in \cite{MR2367118} that under the above assumptions one has $\binom{\mu^+}{\mu}\rightarrow\binom{\mu^+\ \tau}{\mu\quad \mu}$ for every $\tau\in\omega_1$.
Together with the negative result of \cite{GMS} with respect to $\omega_2$, the case of $\omega_1$ seems to be the last case in this context.

We shall prove that $\binom{\mu^+}{\mu}\nrightarrow\binom{\mu^+\ \omega_1}{\mu\quad \mu}$ is consistent for a strong limit singular cardinal $\mu$ with countable cofinality under the assumption that $2^\mu=\mu^+$.
This result can be forced at a limit of measurable cardinals, and also at $\aleph_\omega$.
Thus we obtain a negative answer to \cite[Problem 10]{MR0202613} and we conclude that Jones' result is optimal from \textsf{ZFC} point of view.
We also obtain a negative answer to \cite[Question 6]{GMS} and conclude that in some sense Shelah's result from \cite{MR1606515} is optimal.

Suppose that $\mu\geq\cf(\mu)>\omega$.
As mentioned above, we know that $\binom{\mu^+}{\mu}\nrightarrow\binom{\mu^+\ \omega_1}{\mu\quad \mu}$.
If $c:\mu^+\times\mu\rightarrow 2$ exemplifies this negative relation and we force to make $\cf(\mu)=\omega$ then the coloring $c$ is reinterpreted in the generic extension and it will be polychromatic on old sets of the appropriate size.

The problem is that usually our forcing adds new sets.
For example, if $\mu$ is measurable and we force a Prikry sequence into $\mu$ then new sets of size $\mu$ and $\mu^+$ are added and they might be quite far from old sets of the same size.
In particular, maybe $c$ is monochromatic on these new sets.

One can try to begin with a singular cardinal $\mu$ whose cofinality $\kappa$ is a measurable cardinal, and then to force Prikry into $\kappa$ thus making $\mu$ a singular cardinal with countable cofinality as well.
If one begins with $\mu>\cf(\mu)=\kappa$ and $\kappa$ is measurable then Prikry forcing into $\kappa$ is a bit better since new sets of size $\mu^+$ contain old sets of the same size and this is sufficient for the negative relation.
However, it seems that there is no way to apply a similar argument to sets of size $\mu$.

\begin{proposition}
\label{propoldsets} Suppose that $\mu>\cf(\mu)=\kappa, \kappa$ is a measurable cardinal and $\mathbb{P}$ is Prikry forcing through $\kappa$.
\begin{enumerate}
\item [$(\aleph)$] If $2^\kappa<\theta=\cf(\theta)$ and $A$ is a new set of size $\theta$ then $A$ contains an old subset of size $\theta$.
\item [$(\beth)$] If $\omega<\cf(\theta)\leq\theta<\kappa$ and $A$ is a new set of size $\theta$ then $A$ contains an old subset of size $\theta$.
\item [$(\gimel)$] There exists $A\in[\mu]^\mu$ in the generic extension such that if $B\subseteq A$ and $B$ is an old set then $B$ is bounded in $\mu$.
\end{enumerate}
\end{proposition}

\par\noindent\emph{Proof}. \newline 
Choose a $V$-generic set $G\subseteq\mathbb{P}$.
For the first statement let $A\in V[G]$ be of size $\theta$ and for every $p\in G$ let $A_p=\{\alpha\in{\rm Ord}:p\Vdash\check{\alpha}\in\name{A}\}$.
Notice that $A_p\in V$ whenever one fixes a single condition $p$.
Since $A=\bigcup\{A_p:p\in G\}$ and $|G|\leq 2^\kappa<\theta=\cf(\theta)$ we see that there is a condition $p\in G$ for which $|A_p|=\theta$.
The fact that $A_p$ is forced to be a subset of $A$ concludes the argument.

For the second part suppose that $A$ is a new set of size $\theta$, and without loss of generality $A\subseteq\kappa$, using some one-to-one mapping from $A$ into $\kappa$ and working with the range of this mapping.
Since $\cf(\kappa)=\omega$ in the generic extension and $|A|>\omega$, one can find a subset $B$ of $A$ of the same size which is bounded in $\kappa$, say $\sup(B)=\rho<\kappa$.
Let $p=(s^p,A^p)$ be any condition which forces $B\subseteq A$.
For each $\alpha\in\rho$ let $\varphi_\alpha$ be the formula $\check{\alpha}\in\name{B}$.
For every $\alpha\in\rho$ choose a condition $q_\alpha=(s^p,A^\alpha)$ so that $p\leq^* q_\alpha$ and $q_\alpha$ decides $\varphi_\alpha$.
Let $E=\bigcap\{A^\alpha:\alpha\in\rho\}$, so $E$ belongs to the normal ultrafilter which we use in our forcing.
But now the single condition $(s^p,E)$ determines the elements of $B$ and hence $B\in V$ as required.

Lastly, choose in the ground model an increasing sequence of regular cardinals $(\mu_\delta:\delta\in\kappa)$ such that $\mu=\bigcup_{\delta\in\kappa}\mu_\delta$.
Let $(\rho_i:i\in\omega)$ be a Prikry sequence into $\kappa$.
Set $A=\bigcup\{[\mu_{\rho_i},\mu_{\rho_i}^+):i\in\omega\}$, so $A\in[\mu]^\mu\cap V[G]$.
If $B\subseteq A$ and $B$ is unbounded in $\mu$ then $B\cap[\mu_{\rho_i},\mu_{\rho_i}^+)\neq\varnothing$ for infinitely many $i\in\omega$.
Hence one can recover an infinite subsequence of $(\rho_i:i\in\omega)$ from $B$ by collecting the indices of $\mu_\delta=|\beta|$ for every $\beta\in B$.
By the genericity criterion of Mathias from \cite{MR0332482} we see that $B\notin V$, thus we are done.

\hfill \qedref{propoldsets}

The strategy of using an old coloring is apparently ineffective in the light of the last part of the above proposition.
However, it can be used at non strong limit cardinals.

\begin{claim}
\label{clmnonslimit} Assume that:
\begin{enumerate}
\item [$(a)$] $\mu>\cf(\mu)=\kappa$ and $2^\mu=\mu^+$.
\item [$(b)$] $\kappa$ is measurable and $2^\kappa<\mu$.
\item [$(c)$] There exists some $\theta<\mu$ so that $\kappa<\theta$ and $\binom{\mu^+}{\mu}\nrightarrow\binom{\mu^+\ \omega_1}{\theta\quad \theta}$.
\end{enumerate}
Then one can force $\mu>\cf(\mu)=\omega, 2^\mu=\mu^+$ and $\binom{\mu^+}{\mu}\nrightarrow\binom{\mu^+\ \omega_1}{\theta\quad \theta}$ so one obtains $\binom{\mu^+}{\mu}\nrightarrow\binom{\mu^+\ \omega_1}{\mu\quad \mu}$ as well.
\end{claim}

\par\noindent\emph{Proof}. \newline 
Let $\mathbb{P}$ be Prikry forcing through $\kappa$ and let $G\subseteq\mathbb{P}$ be $V$-generic.
Fix, in the ground model, a coloring $c:\mu^+\times\mu\rightarrow 2$ which exemplifies the negative relation $\binom{\mu^+}{\mu}\nrightarrow\binom{\mu^+\ \omega_1}{\theta\quad \theta}$.
By abuse of notation let us denote $c_G$ by $c$.

Assume that $A\in[\mu^+]^{\mu^+}$ and $B\in[\mu]^\theta$ in $V[G]$.
By Proposition \ref{propoldsets} one can find $a\in[A]^{\mu^+}\cap V$ and $b\in[B]^\theta\cap V$.
By $(c)$ we know that $1\in c''(a\times b)$ and then $1\in c''(A\times B)$.
Similarly, if $A\in[\mu^+]^{\omega_1}$ and $B\in[\mu]^\theta$ then one can choose $a\in[A]^{\omega_1}$ and $b\in[B]^\theta$ from the ground model.
It follows that $0\in c''(a\times b)\subseteq c''(A\times B)$, so we are done.

\hfill \qedref{clmnonslimit}

Assumptions in the spirit of the above claim can be forced by adding many Cohen sets to some relatively small cardinal.
This idea is applicable to non strong limit cardinals, and an interesting example will be proved anon.
We mention here \cite[Problem 14]{MR0202613}, which asks whether $\binom{\aleph_{\omega+1}}{\aleph_{\omega+1}}\rightarrow\binom{\aleph_{\omega+1}\ \omega_1}{\aleph_\omega\quad \aleph_\omega}$.

\begin{claim}
\label{clm14} It is consistent that $\mu>\cf(\mu)=\omega, 2^\omega<\mu$, but $\binom{\mu^+}{\mu^+}\nrightarrow\binom{\mu^+\ \omega_1}{\mu\quad \mu}$ and even $\binom{\mu^+}{\mu^+}\nrightarrow\binom{\mu^+\ \omega}{\mu\quad \mu}$.
\end{claim}

\par\noindent\emph{Proof}. \newline 
Begin with $\mu>\cf(\mu)=\kappa$ such that $\kappa$ is measurable and $2^\kappa=\kappa^+$.
Choose $\theta=\cf(\theta)$ such that $\kappa^+<\theta<\mu$.
Let $\mathbb{Q}$ be $Add(\theta,\mu^+)$ and let $G\subseteq\mathbb{Q}$ be $V$-generic.
By \cite[Theorem 7.4]{MR401472} we have $\binom{\mu^+}{\mu^+}\nrightarrow\binom{\theta^+}{\theta}$ in $V[G]$.

Working in $V[G]$, let $\mathbb{P}$ be Prikry forcing through $\kappa$ (notice that $\kappa$ remains measurable in $V[G]$).
Let $H\subseteq\mathbb{P}$ be $V[G]$-generic.
The argument in the proof of Claim \ref{clmnonslimit} gives $\binom{\mu^+}{\mu^+}\nrightarrow\binom{\mu^+\ \theta^+}{\theta\quad \theta}$ in $V[G][H]$.

Moreover, suppose that $\theta<\lambda<\mu$ and $\cf(\lambda)=\omega_1$.
By monotonicity we have $\binom{\mu^+}{\mu^+}\nrightarrow\binom{\mu^+\ \lambda}{\theta\quad \theta}$ in $V[G][H]$, and since $2^{\omega_1}<\mu$ we see that $\binom{\mu^+}{\mu^+}\nrightarrow\binom{\mu^+\ \omega_1}{\theta\quad \theta}$ as well.
Another application of monotonicity yields $\binom{\mu^+}{\mu^+}\nrightarrow\binom{\mu^+\ \omega_1}{\mu\quad \mu}$ in $V[G][H]$, as desired.
The same argument gives the negative relation $\binom{\mu^+}{\mu^+}\nrightarrow\binom{\mu^+\ \omega}{\mu\quad \mu}$ if one chooses $\lambda>\theta$ with $\cf(\lambda)=\omega$, so the proof is accomplished.

\hfill \qedref{clm14}

Observe that $2^{\cf(\mu)}<\mu$ is satisfied in the generic extension, and for every specific $\chi<\mu$ one can force as in the above claim while retaining $2^\chi<\mu$ by choosing $\theta\in(\chi,\mu)$.
However, the above idea seems to be inapplicable if one considers strong limit cardinals.
In particular, the following claim shows that negative assumptions in the ground model with respect to some $\theta<\mu$ as needed for the above proof are impossible.

\begin{claim}
\label{clmslimit} Suppose that $\mu$ is a strong limit cardinal and $\theta<\mu$. Then $\binom{\mu^+}{\mu}\rightarrow\binom{\mu^+\ \mu^+}{\mu\quad \theta}$.
\end{claim}

\par\noindent\emph{Proof}. \newline 
Suppose that $c:\mu^+\times\mu\rightarrow 2$.
Let $\kappa=\cf(\mu)\leq\mu$ and let $(\mu_\gamma:\gamma\in\kappa)$ be an increasing sequence of cardinals such that $\mu=\bigcup_{\gamma\in\kappa}\mu_\gamma$.
For every $\alpha\in\mu^+$ and each $i\in\{0,1\}$ define:
$$
A_{\alpha i}=\{\beta\in\mu:c(\alpha,\beta)=i\}.
$$
Let $\mathcal{F}=\{A_{\alpha 1}:\alpha\in\mu^+\}$.
By \cite[Lemma 1]{MR2367118} either there is some $B\in[\mu^+]^{\mu^+}$ such that $|\bigcap\{\mu-A_{\alpha 1}:\alpha\in B\}|=\mu$ or there exists $A\in[\mu^+]^{\mu^+}$ such that $\{A_{\alpha 1}:\alpha\in A\}$ is a $\mu$-uniform filter base.

Lest the first option obtains we have a $0$-monochromatic product $B\times C$ of size $\mu^+\times\mu$, stipulating $C=\bigcap\{\mu-A_{\alpha 1}:\alpha\in B\}$.
If the second option holds, choose for every $\alpha\in A$ an ordinal $\gamma(\alpha)\in\kappa$ such that $|A_{\alpha 1}\cap\mu_{\gamma(\alpha)}|\geq\theta$.
Since $|A|=\mu^+$ we may assume, without loss of generality, that $\gamma(\alpha)=\gamma$ for every $\alpha\in A$ and some fixed $\gamma\in\kappa$.
Since $\mu_\gamma<\mu$ and $\mu$ is strong limit we see that $|\mathcal{P}(\mu_\gamma)|<\mu<\mu^+$.
Hence without loss of generality there is a fixed $a\subseteq\mu_\gamma$ such that $|a|\geq\theta$ and $A_{\alpha 1}\cap\mu_\gamma=a$ for every $\alpha\in A$.
Verify that $A\times a$ is $1$-monochromatic and conclude that $\binom{\mu^+}{\mu}\rightarrow\binom{\mu^+\ \mu^+}{\mu\quad \theta}$ as desired.

\hfill \qedref{clmslimit}

Although we cannot use old colorings in the context of strong limit cardinals, we can still exploit the density arguments for sets of size $\omega_1$.
The following is the main result of this section:

\begin{theorem}
\label{thmmt} It is consistent that $\mu$ is a strong limit singular cardinal of countable cofinality, $2^\mu=\mu^+$ and $\binom{\mu^+}{\mu}\nrightarrow\binom{\mu^+\ \omega_1}{\mu\quad \mu}$.
\end{theorem}

\par\noindent\emph{Proof}. \newline 
Let $\mu$ be a measurable cardinal, assume that $2^\mu=\mu^+$ and let $\mathbb{P}$ be Prikry forcing into $\mu$.
We claim that $\binom{\mu^+}{\mu}\nrightarrow\binom{\mu^+\ \omega_1}{\mu\quad \mu}$ in the generic extension by $\mathbb{P}$.
Choose a generic subset $G\subseteq\mathbb{P}$.
We shall define a coloring $c:\mu^+\times\mu\rightarrow 2$ in $V[G]$.

As a first step we enumerate $[\mu^+]^\mu$ by $\{B_\alpha:\alpha\in\mu^+\}$.
Now for each $\alpha\in(0,\mu^+)$ we enumerate $\{B_\beta:\beta\in\alpha\}$ in such a way that the order-type will be $\mu$, say $\{B_{\alpha\varepsilon}:\varepsilon\in\mu\}$.
If $\alpha\in\mu$ then we use repetitions.
Similarly, for each $\alpha\in(0,\mu^+)$ we list the ordinals of $\alpha$ by an enumeration of order-type $\mu$, say $\{\alpha_\eta:\eta\in\mu\}$.
We emphasize that this enumeration of the ordinals of $\mu^+$ is done already in the ground model.

As a second step we choose for every $\alpha\in\mu^+$ and each $B_{\alpha\varepsilon}$ an ordinal $\beta\in B_{\alpha\varepsilon}$ with the goal of setting $c(\alpha,\beta)=1$.
We do this, however, with some control.
Fix $\alpha\in\mu^+$ and choose for each $\varepsilon\in\mu$ an ordinal $\gamma_{\alpha\varepsilon}$ so that $\gamma_{\alpha\varepsilon}\in B_{\alpha\varepsilon}-\{\gamma_{\alpha_\eta\zeta}: \eta<\varepsilon,\zeta\leq\varepsilon\}$.
Notice that we remove a small set from $B_{\alpha\varepsilon}$ and $|B_{\alpha\varepsilon}|=\mu$, so the choice is possible. 
Define:
$$
c(\alpha,\beta)=1 \Leftrightarrow \exists\varepsilon\in\mu, \beta=\gamma_{\alpha\varepsilon}.
$$
We must show that $c$ has no $0$-monochromatic product of size $\mu^+\times\mu$ and no $1$-monochromatic product of size $\omega_1\times\mu$.

For the first mission suppose that $A\in[\mu^+]^{\mu^+}$ and $B\in[\mu]^\mu$.
Let $\alpha\in A$ be such that $B\in\{B_\beta:\beta\in\alpha\}$ and let $\varepsilon\in\mu$ be such that $B=B_{\alpha\varepsilon}$.
Let $\beta=\gamma_{\alpha\varepsilon}$, so $\beta\in B_{\alpha\varepsilon}=B$.
By definition, $c(\alpha,\beta)=1$ and hence $A\times B$ is not $0$-monochromatic.

For the second mission assume that $A\in[\mu^+]^{\omega_1}$ and $B\in[\mu]^\mu$.
Assume toward contradiction that $A\times B$ is $1$-monochromatic.
Choose $A'\subseteq A$ such that $|A'|=\aleph_1$ and $A'\in V$.
Define $S=\{\eta\in\mu:\exists\gamma,\delta\in A',\gamma=\delta_\eta\}$.
Notice that $S\subseteq\mu$ and $S\in V$ since $A'\in V$.
In particular, $S$ is bounded in $\mu$, so let $\rho=\sup(S)<\mu$.
Define $T=\{\beta\in\mu:\forall\alpha\in A',c(\alpha,\beta)=1\}$.
Observe that $T\supseteq B$ since $A'\subseteq A$ and for every $\alpha\in A,\beta\in B$ we have $c(\alpha,\beta)=1$.
We shall prove that $|T|<\mu$ thus arriving at a contradiction since $|B|=\mu$.

Firstly we observe that if $\alpha\in\mu^+$ and $\eta<\varepsilon<\mu$ and $\gamma(\alpha_\eta,\zeta)=\gamma(\alpha,\varepsilon)$ then $\varepsilon<\zeta$.
To see this recall that $\alpha_\eta<\alpha$, so if $\zeta\leq\varepsilon$ then we required in the choice of the $\gamma$s that $\gamma(\alpha_\eta,\zeta)\neq\gamma(\alpha,\varepsilon)$.
Secondly, choose some $\beta\in T$ and recall that $c(\alpha,\beta)=1$ for every $\alpha\in A'$.
It follows that for each $\alpha\in A'$ there exists a unique ordinal $\varepsilon(\alpha)\in\mu$ such that $\beta=\gamma_{\alpha\varepsilon(\alpha)}$.
We claim that there must be some $\alpha\in A'$ for which $\varepsilon(\alpha)\in\rho$.
Suppose not, and choose $\sigma,\delta\in A'$ such that $\sigma<\delta$.
Let $\eta\in\mu$ be so that $\sigma=\delta_\eta$.
By the definition of $\rho$ we see that $\eta<\rho$, so $\eta<\varepsilon(\delta)$ by our assumption.
Since both $\gamma_{\sigma\varepsilon(\sigma)}=\beta$ and $\gamma_{\delta\varepsilon(\delta)}=\beta$ we see that $\gamma_{\sigma\varepsilon(\sigma)}=\gamma_{\delta\varepsilon(\delta)}$ and by the above observation we have $\varepsilon(\sigma)>\varepsilon(\delta)$.
Since $A'$ is infinite, if we choose an increasing sequence $\langle\sigma_n:n\in\omega\rangle$ of elements of $A'$ we produce an infinite decreasing sequence of ordinals $\langle\varepsilon(\sigma_n):n\in\omega\rangle$, which is an absurd.

Therefore, for each $\beta\in T$ we choose $\alpha\in A'$ such that $\beta=\gamma_{\alpha\varepsilon(\alpha)}$ and $\varepsilon(\alpha)<\rho$.
It follows now from the definition of the set $T$ that $T\subseteq\{\gamma_{\alpha\varepsilon}:\alpha\in A', \varepsilon<\rho\}$ and hence $|T|\leq|\rho|\cdot\aleph_1<\mu$, so we are done.

\hfill \qedref{thmmt}

The above theorem was proved using Prikry forcing, and this is probably the simplest way to carry out the argument.
But the method itself is a bit more general, and in particular applies to a wider collection of forcing notions.
If $\mathbb{P}$ forces countable cofinality to $\mu$ where $\mu\geq\cf(\mu)>\omega$ in the ground model and every new set of size $\aleph_1$ contains an old set of size $\aleph_1$ then a similar proof works.
Consequently, one can use other Prikry-type forcing notions.
An answer to \cite[Problems 10 and 14]{MR0202613} can be given now using Magidor's method from \cite{MR491183} and \cite{MR491184} to singularize a measurable cardinal with interleaved collapses.

\begin{corollary}
\label{coralephomega} It is consistent that $\binom{\mu^+}{\mu}\nrightarrow\binom{\mu^+\ \omega_1}{\mu\quad \mu}$ where $\mu=\aleph_\omega$ is strong limit and $2^\mu=\mu^+$.
Similarly, $\binom{\mu^+}{\mu^+}\nrightarrow\binom{\mu^+\ \omega_1}{\mu\quad \mu}$ is consistent at $\mu=\aleph_\omega$ under the same assumptions.
\end{corollary}

\hfill \qedref{coralephomega}

The following remark which sheds some light on the role of $\omega_1$ in the above statements.
The general idea behind the proof of the main result it to begin with a negative relation at $\mu$, to force a desired property of $\mu$ (here it is countable cofinality) in such a way that the negative relation is preserved.
Thus we commence with a measurable cardinal $\mu$ with $2^\mu=\mu^+$ so $\binom{\mu^+}{\mu}\nrightarrow\binom{\mu^+\ \omega_1}{\mu\quad \mu}$ in the ground model, and we force Prikry in order to singularize $\mu$ while keeping the negative relation.

But under the above assumption we know that $\binom{\mu^+}{\mu}\nrightarrow\binom{\mu^+\ \omega}{\mu\quad \mu}$ as well in the ground model.
When we force with Prikry forcing we secure the negative relation $\binom{\mu^+}{\mu}\nrightarrow\binom{\mu^+\ \omega_1}{\mu\quad \mu}$ but we also get the positive relation $\binom{\mu^+}{\mu}\rightarrow\binom{\mu^+\ \omega}{\mu\quad \mu}$ in the generic extension.
One may wonder what is the difference between $\omega$ and $\omega_1$ in this context, and at least one aspect of the answer becomes clearer.
New sets of size $\aleph_1$ in the Prikry extension contain old sets of the same size.
Contrariwise, new $\omega$-sequences cannot be approximated by infinite sets from the ground model.
This fact explains the peculiarity of singular cardinals with countable cofinality with respect to the unbalanced relation discussed in this paper.

We conclude with a couple of open problems.
The first one is about the possibility of the positive direction.

\begin{question}
\label{qpositive} Is it consistent that the positive relation $\binom{\mu^+}{\mu}\rightarrow\binom{\mu^+\ \omega_1}{\mu\quad \mu}$ holds under \textsf{GCH} where $\mu>\cf(\mu)=\omega$?
\end{question}

For the second question let us indicate that large cardinals are not necessary in order to force $\binom{\mu^+}{\mu}\nrightarrow\binom{\mu^+\ \omega_1}{\mu\quad \mu}$ under the above assumptions.
For example, if one forces with Namba forcing over the constructible universe then every singular cardinal $\mu$ so that $\cf^V(\mu)=\omega_2$ will satisfy $\binom{\mu^+}{\mu}\nrightarrow\binom{\mu^+\ \omega_1}{\mu\quad \mu}$ in the generic extension.
However, by assuming the existence of large cardinals one can force a universe in which every singular strong limit cardinal with countable cofinality is ex-measurable, and then $\binom{\mu^+}{\mu}\nrightarrow\binom{\mu^+\ \omega_1}{\mu\quad \mu}$ holds globally at every such cardinal.
We do not know whether large cardinals are indispensable for this result:

\begin{question}
\label{qglobal} Can one force without large cardinals that \textsf{GCH} holds, and $\binom{\mu^+}{\mu}\nrightarrow\binom{\mu^+\ \omega_1}{\mu\quad \mu}$ whenever $\mu>\cf(\mu)=\omega$?
\end{question}

\newpage 

\bibliographystyle{amsplain}
\bibliography{arlist}

\end{document}